\documentstyle[12pt]{article}

\newtheorem{lemma}{Lemma}
[section]

\newtheorem{thm}{Theorem}
[section]
\newtheorem{cor}{Corollary}
[section]

\newcommand{\mysection}[1]{\section{#1}\setcounter{equation}{0}}

\newcommand{\la}{\lambda}

\newcommand{\ptl}{\partial}

\newcommand{\ra}{\rightarrow}
\newcommand{\lra}{\longrightarrow}

\def\a{\alpha}

\def\R2{{\Bbb R}^2}
\def\Bbb{\mathbb}

\def\beq{\begin{equation}}
\def\eeq{\end{equation}}
\def\ba{\begin{array}}
\def\ea{\end{array}}
\def\barr{\begin{eqnarray}}
\def\earr{\end{eqnarray}}

\def\part{\partial}
\def\fr{\frac}
\def\ov{\overline}

\newcommand{\ls}{\setlength{\baselineskip}{18pt}
                      \setlength{\parskip}{3mm} }

\title {The $C^{\a}$ regularity of a class of ultraparabolic equations}
\author{ ZHANG Liqun
\thanks{The author currently is working at NSFC .}}

\date{Institute of Mathematics, AMSS, Academia Sinica,
Beijing}

\begin{document}

\maketitle

\begin{abstract}
We prove the $C^{\a}$ regularity for weak solutions to a class of
ultraparabolic equation, with measurable coefficients. The results
generalized our recent $C^{\a}$ regularity results of Prandtl's
system to high dimensional cases.

\end{abstract}

{\small keywords: Ultraparabolic equations, Moser iteration,
$C^{\a}$ regularity }

\pagenumbering{arabic}
\mysection{Introduction} \ls \noindent

The ultraparabolic equation arises in many applications, for
example, fluid dynamics, mathematical finance, degenerated
diffusion process, etc. There are more and more studies on this
problem in recent years. The regularity of this type of equation
becomes interesting since it has some special algebraic structures
and is degenerated. It is still unclear in general, whether the
interior $C^{\a}$ regularity results hold for weak solutions of
the ultraparabolic equations with bounded measurable coefficients
like the parabolic cases.

In the study of boundary layer problem, we obtained the existence
global weak solution in the class that Oleinik considered under
the assumption that the pressure is favorable [13]. One of the
interesting question is whether the weak solution is actually
smooth?  In particular, in the two dimensional Prandtl's system
with constant pressure, in the Crocco variable, we can deduce the
following equation
$$
\displaystyle \fr{\partial \,u}{\partial \,t}+y \fr{\partial
\,u}{\partial \,x}-u^2\fr{\partial^2 \,u}{\partial
\,y^2}=0.\leqno(1.1)
$$
This is of strong degenerated parabolic type equations, more
precisely, an ultraparabolic type equation. However, it satisfies
the well known H\"ormander's hypoelliticity conditions, which
sheds lights on the smoothness of weak solutions. It is proved by
Polidoro and  Ragusa [12] that the weak solution is in the
$C^{\a}$ class, if the coefficient is in the class of VMO. That
is, if the weak solution $u$ is also in the class of VMO, then the
solution $u$ is $C^{\a}$ continuous and then the result of Xu
[15], or Bramanti, Cerutti and Manfredini [1] implies that the
solution is also $C^{\infty}$ smooth. In our case, however, we can
only prove that the weak solution is in the class of BV space. It
is interesting if the weak solution of equation (1.1) is still
smooth when the coefficient is only measurable functions.

On the other hand, equation (1.1) has the divergent form if we
replace $u$ by $\fr 1 u$. A recent paper by Pascucci and Polidoro
[11] proved that the Moser iterative method still works for a lass
of ultraparabolic equations with measurable coefficients including
equation (1.1). Their results showed that for a non-negative
sub-solution $u$ of (1.1), the $L^{\infty}$ norm of $u$ is bounded
by the $L^p$ norm of $u$ ($p \ge 1$). This is a very important
step toward the final solution of regularity of the ultraparabolic
equations.

We proved in [14] that the weak solution that we obtained in [13]
of (1.1) is of $C^{\a}$ class and then $u$ is smooth. In this
paper, we are concerned with the $C^{\a}$ regularity of solutions
of the ultraparabolic equations. We shall generalize the result in
[13] to high dimensional cases in this paper.

We consider a class of Komogorov-Fokker-Planck type operator on
${R}^{N+1}$.
$$
\displaystyle  Lu \equiv
\sum_{i,j=1}^{m_0}{\ptl_{x_j}\,(a_{ij}(t,x)\ptl_{x_i}\,u
)}+\sum_{i,j=1}^N b_{ij}x_i {\ptl_{x_j}\, u }- {\ptl_t
\,u}=0,\leqno(1.2)
$$
where $(x,t)\in { R}^{N+1}$, $1\leq m_0 \leq N$, and $b_{ij}$ is
constant for every $i,j= 1,\cdots, N$. We make the following
assumptions on the coefficients of $L$:

$(H_1)$  $a_{ij}=a_{ji} \in L^{\infty} ({R}^{N+1})$ and there
exists a $\lambda >0$ such that
$$
\fr{1}{\lambda}\sum_{i=1}^{m_0}\xi_i^2 \leq \sum_{i,j=1}^{m_0}
a_{ij}(t,x)\xi_i \xi_j \leq {\lambda}\sum_{i=1}^{m_0}\xi_i^2
$$
for every $(t,x)\in {R}^{N+1}$, and $\xi \in {R}^{m_0}$.

$(H_2)$ The matrix $B=(b_{ij})_{N \times N}$ has the form
$$\left(
\begin{array}{cccccc}
0 & {B_1} & 0 & \cdots & 0 \\
0 & 0  & {B_2} & \cdots & 0 \\
\vdots & \vdots & \vdots & \ddots &\vdots \\
0 & 0 & 0 & \cdots & {B_d} \\
0 & 0 & 0 & \cdots & 0
\end{array}
\right)
$$
where $B_k$ is a matrix $m_{k-1}\times m_{k}$ with rank $m_k$ and
$m_0\geq m_1\geq \cdots\geq m_d$, $m_0+m_1+\cdots+m_d=N$.
$||B||\leq \la$ where the norm $||\cdot||$ is in the sense of
matrix norm.

The requirements of matrix $B$ in $(H_2)$ ensures that the
operator $L$ with constant $a_{ij}$ satisfies the well-known
H\"omander's hypoellipticity condition.

The Schauder type estimate of (1.2) has been obtained. Besides, the
regularity of weak solutions have been studied by Bramanti, Cerutti
and Manfredini [1], Manfredini and Polidoro [7], Polidoro and Ragusa
[12] assuming a weak continuity on the coefficient $a_{ij}$. It is
quite interesting whether the weak solution has H\"older regularity
under the assumption $(H_1)$ on $a_{ij}$. The first advances is the
work of A. Pascucci and S. Polidoro [11] who proved that Moser
iteration method still works for equation (1.2). One of the approach
to the H\"older estimates is to obtain the Harnack type inequality.
In the case of elliptic equation with measurable coefficients, the
Harnack inequality is obtained by J. Moser [8] via an estimate of
BMO functions due to F. John and L. Nirenberg together with the
Moser iteration method. J. Moser [9] also obtained the Harnack
inequality for parabolic equations with measurable coefficients by
generalizing the John-Nirenberg estimates to the parabolic case.
Another approach to the H\"older estimates is given by S. N.
Kruzhkov [6], [7] based on the Moser iteration to obtain a local
priori estimates, which provides a short proof for the parabolic
equations.

We prove a Poincare type inequality for non-negative weak
sub-solutions of (1.2). Then we apply it to obtain a local priori
estimates which implies the H\"older estimates for ultraparabolic
equation (1.2).

Let $D_{m_0}$ be the gradient with respect to the variables $x_1,
x_2,\cdots, x_{m_0}$. And
$$Y=\sum_{i,j=1}^N b_{ij}x_i {\ptl_{x_j}}- {\ptl_t}.$$
We say that $u$ is a weak solution of (1.2) if it satisfies (1.2)
in the distribution sense and $u$, $D_{m_0}u$, $Yu \in L^2_{loc}$.

Our main result is the following theorem.

\begin{thm}
Under the assumptions $(H_1)$ and $(H_2)$, the weak solution of
(1.2) is H\"older continuous.
\end{thm}

\mysection{Some Preliminary Results}

One of the important feature of equation (1.2) is that the
fundamental solution can be written down explicitly if the
coefficients $a_{ij}$ is constant, (see [1], [4]). Besides, there
are some geometric and algebraic structures in the space $R^{N+1}$
induced by the constant matrix $B$ (see for instance, [1]).

Let $E(\tau)=exp(-\tau B^T)$, where $E(\tau)$ is a polynomial of
degree $d$ in $\tau$ with $N \times N$ matrices coefficients. For
$(t,x), (\tau,y) \in R^{N+1}$, set $$(t,x)\circ
(\tau,y)=(t+\tau,y+E(\tau)x).$$ Then $(R^{N+1}, \circ)$ is a group
with neutral element $(0,0)$; the inverse of an element $(t,x)$ is
$(t,x)^{-1}=(-t,-E(-t)x)$. The left translation by $(\tau,y)$
given by
$$(t,x)\mapsto (\tau,y)\circ (t,x),$$
is a invariant translation to operator $L$ when coefficient
$a_{ij}$ is constant.

The associated dilation to operator $L$ with constant coefficient
$a_{ij}$ is given by
$$
\delta_{\la}=diag(\la^2,\la I_{m_0},\la^3
I_{m_1},\cdots,\la^{2d+1}I_{m_d}),
$$
where $I_{m_k}$ denotes the $m_k\times m_k$ identity matrix. Then
the operator is homogeneous of degree 2 with respect to the
dilation $\delta_{\la}$. Let
$$Q=m_0+3m_1+\cdots+(2d+1)m_d.$$
Then the number $Q+2$ is usually called the homogeneous dimension
of $R^{N+1}$ with respect to the dilation $\delta_{\la}$.

The norm in $R^{N+1}$, related to the group of translations and
dilation to the equation is defined by $$||(t,x)||=r$$ if $r$ is the
unique positive solution to the equation
$$
\fr{x_1^2}{r^{2\a_1}}+\fr{x_2^2}{r^{2\a_2}}+\cdots+\fr{x_N^2}{r^{2\a_N}}
+\fr{t^2}{r^4}=1,
$$
where $(t,x) \in R^{N+1}\setminus \{0\}$ and
$$
\a_1=\cdots=\a_{m_0}=1, \quad
\a_{m_0+1}=\cdots=\a_{m_0+m_1}=3,\cdots,
$$
$$
\a_{m_0+\cdots+m_{d-1}+1}=\cdots=\a_N=2d+1.
$$
And $||(0,0)||=0$. The balls at a point $(t_0,x_0)$ is defined by
$${\cal B}_r(t_0,x_0)=\{(t,x)|\quad ||(t,x)^{-1}\circ (t_0,x_0)||\leq r\}.$$
Let
$${\cal B}^-_r(t_0,x_0)={\cal B}_r(t_0,x_0)\cap\{t<t_0\}.$$

For convenience, we sometimes use the cube replace the balls. The
cube at point $(0,0)$ is given by
$$
{\cal C}_r(0,0)=\{(t,x)|\quad |t|\leq r^2,\quad |x_1|\leq
r^{\a_1}, \cdots, |x_N|\leq r^{\a_N}\}.
$$
It is easy to see that there exists a constant $\Lambda$ such that
$$
{\cal C}_{\fr r \Lambda}(0,0)\subset{\cal B}_r(0,0)\subset{\cal
C}_{\Lambda r}(0,0),
$$
where $\Lambda$ only depends on $B$ and $N$.

When the matrix $(a_{ij})_{N \times N}$ is of constant matrix, we
denoted it by $A_0$. Then the operator $L$ takes the form
$$
L_0=div(A_0D)+Y.
$$
We let $z=(t,x)$. The fundamental solution $\Gamma_0(\cdot,\zeta)$
of $L_0$ with pole in $\zeta\in R^{N+1}$ has been constructed (see
[1]) as follows:
$$
\Gamma_0(z,\zeta)=\Gamma_0(\zeta^{-1}\circ z,0), \qquad z, \zeta
\in R^{N+1},\quad z \neq \zeta,
$$
And $\Gamma_0(z,0)$ can be written down explicitly. There are some
basic estimates for $\Gamma_0$
$$
\Gamma_0(z,\zeta)\leq c ||\zeta^{-1}\circ z||^{-Q},
$$
$$
|\ptl_{x_i}\,\Gamma_0(z,\zeta)|\leq c ||\zeta^{-1}\circ
z||^{-Q-1},
$$
where $i=1,\cdots,m_0$.

A weak sub-solution of (1.2) in a domain $\Omega$ is a function
$u$ such that $u$, $D_{m_0}u$, $Yu \in L^2_{loc}(\Omega)$ and for
any $\phi \in C^{\infty}_0(\Omega)$, $\phi \geq 0$,
$$
\int_{\Omega} \phi Yu-(Du)^T AD\phi \geq 0. \leqno(2.1)
$$

A result of Pascucci and Polidoro obtained by using the Moser's
iterative method (see [11]) states as follows.

\begin{lemma}
Let $u$ be a non-negative weak sub-solution of (1.2) in $\Omega$.
Let $(t_0,x_0)\in \Omega$ and $\overline{{\cal
B}^-_r(t_0,x_0)}\subset \Omega$ and let $p \geq 1$. Then there
exists a positive constant $c$ which depends only on $\la$ and the
homogeneous dimension $Q$ such that, for $0 < r\leq 1$
$$
\sup_{{\cal B}^-_{\fr r 2}(t_0,x_0)} u^p \leq \fr
{c}{r^{Q+2}}\int_{{\cal B}^-_r(t_0,x_0)} u^p,\leqno(2.2)
$$
provided that the last integral converges.
\end{lemma}

We copy a classical potential estimates (see [3]) here to prove
the Poincare type inequality.

\begin{lemma}
Let $\a \in (0, Q+2)$ and $G \in C(R^{N+1}\setminus \{0\})$ be a
$\delta_{\la}$-homogeneous function of degree $\a-Q-1$. If $f \in
L^p(R^{N+1})$ for some $p \in (0,\infty)$, then
$$
G_f(z)\equiv \int_{R^{N+1}} G(\zeta ^{-1}\cdot z)f(\zeta)d\zeta,
$$
is defined almost everywhere and there exists a constant
$C=C(Q,p)$ such that
$$
||G_f||_{L^q(R^{N+1})}\leq C \max_{||z||=1} |G(z)|\quad
||f||_{L^p(R^{N+1})},\leqno(2.3)
$$
where $q$ is defined by
$$
\fr 1q =\fr 1p-\fr{\a}{Q+2}.
$$
\end{lemma}

\mysection{Proof of Main Theorem}

We want to obtain a local estimates of solutions of the equation
(1.2), for instant, at point $(t_0,x_0)$. Since the equation (1.2)
is invariant under the left group translation when $a_{ij}$ is
constant, we may consider the estimates at a ball centered at
$(0,0)$. We mainly prove the following Lemma 3.4 which is
essential in the oscillation estimates in Kruzhkov's approaches in
parabolic case. Then the $C^{\a}$ regularity result follows easily
by the standard arguments.

For convenience, in the following discussion, we let
$x'=(x_1,\cdots,x_{m_0})$ and $x=(x', \overline x)$. We consider the
estimates in the following cube, instead of ${\cal B}^-_r$,
$$
{\cal C}_r=\{(t,x)| -r^2\leq t \leq 0,\quad |x'|\leq r,
|x_{m_0+1}|\leq \la N^2 r^{3} \cdots, |x_N|\leq \la N^2
r^{2d+1}\}.
$$
Let
$$
S_r=\{ \overline x\quad|\quad|x_{m_0+1}|\leq \la N^2 r^{3} \cdots,
|x_N|\leq \la N^2 r^{2d+1}\}.
$$

Let $0<\a, \beta<1$ be constant and
$$
K_r=\{x'|\quad |x'|\leq r \},
$$
$$
S_{\a r}=\{ \overline x\quad|\quad|x_{m_0+1}|\leq \la N^2 {(\a
r)}^{3} \cdots, |x_N|\leq \la N^2 {(\a r)}^{2d+1}\},
$$
$$
K_{\a r}=\{x'|\quad |x'|\leq \a r \}.
$$
Now for fixed $t$, let
$$
{\cal N}_t=\{(x',\overline x)\in K_{\beta r}\times S_{\beta
r},\quad u \geq h\}.
$$

In the following discussions, we sometimes abuse the notations of
${\cal B}^-_r$ and ${\cal C}_r$, since there are equivalent.

\begin{lemma}
Suppose that $u(t,x)\geq 0$ be a solution of equation (1.2) in
${\cal B}^-_r$ centered at $(0,0)$ and
$$
mes\{(t,x)\in {\cal B}^-_r, \quad u \geq 1\} \geq \fr 1 2 mes
({\cal B}^-_r).
$$
Then there exist constant $\a$, $\beta$ and $h$, $0<\a, \beta,
h<1$ which only depend on $\la$ and $N$ such that for all $t\in
(-\a r^2,0)$,
$$
mes\{{\cal N}_t\} \geq \fr {1}{11}mes\{ K_{\beta r}\times S_{\beta
r}\}.
$$
\end{lemma}
{\it Proof:} Let
$$
v=\ln^+(\fr{1}{u+h^{\fr 9 8}}),
$$
where $h$ is a constant $0<h<1$ to be determined later. Then $v$
at points where $v$ is positive, satisfies
$$
\displaystyle
\sum_{i,j=1}^{m_0}{\ptl_{x_j}\,(a_{ij}(t,x)\ptl_{x_i}\,v
)}-(Dv)^TA Dv+x^T B Dv - {\ptl_t \,v}=0.\leqno(3.1)
$$
Let $\eta(x')$ be a smooth cut-off function so that
$$
\eta(x')=1,\quad \hbox {for} \quad |x'|< \beta r;
$$
$$
\eta(x')=0,\quad \hbox {for} \quad |x'|\geq r.
$$
Moreover, $0\leq\eta \leq 1$ and $|D_{m_0} \eta|\leq \fr
{2m_0}{(1-\beta)r}$.

Multiplying $\eta^2(x')$ to (3.1) and integrating by parts
$$
\ba{lllll} \int_{K_{\beta r}}\int_{S_{\beta r}} v(t,x',\overline
x)d \overline x dx' +\fr {1}{2\la}\int_\tau^t \int_{K_{\beta
r}}\int_{S_{\beta r}}
|D_{m_0}v|^2d \overline x dx'dt \\ \\
\leq \fr {C} {\beta^{Q}(1-\beta)^2}mes(S_{\beta r})mes(K_{\beta
r})+\int_\tau^t
\int_{K_{r}}\int_{S_{\beta r}} \eta^2 x^TBDv d \overline xdx' dt \\ \\
\quad +\int_{K_{r}}\int_{S_{\beta r}} v(\tau,x',\overline x)d
\overline x dx', \ea\leqno(3.2)
$$
where $C$ only depends on $\la$ and $N$.

Integrating by parts, we have for any $i,j$
$$
\int_{K_{r}}\int_{S_{\beta r}} \eta^2 x_ib_{ij}\ptl_{x_j}v d
\overline x dx'\leq \fr {r^{-2}} {4N^2} {\beta}^{-2Q} \ln(\fr
{1}{h^{\fr 98}}) mes(S_{\beta r})mes(K_{\beta r}).\leqno(3.3)
$$
Then
$$
\int_\tau^t\int_{K_{r}}\int_{S_{\beta r}} \eta^2 x^TBDv d
\overline x dx'dt\leq \fr {1}{4} {\beta}^{-2Q} \ln(\fr {1}{h^{\fr
98}}) mes(S_{\beta r})mes(K_{\beta r}).\leqno(3.4)
$$
We shall estimate the measure of the set ${\cal N}_t$. Let
$$
\mu(t)=mes\{(x',\overline x)|\quad x'\in K_r,\quad \overline x \in
S_{r}, \quad u\geq 1\}.
$$
By our assumption, (for convenience, we may let ${\cal B}^-_r$ be
replaced by ${\cal C}^-_r$), for $0<\a< \fr 12$
$$
\fr 12 r^2 mes(S_{r})mes(K_{r})\leq \int_{-r^2}^0
\mu(t)dt=\int_{-r^2}^{-\a r^2}\mu(t)dt+\int_{-\a r^2}^{0}\mu(t)dt.
$$
That is
$$
\int_{-r^2}^{-\a r^2}\mu(t)dt\geq (\fr 12-\a)r^2
mes(S_{r})mes(K_{r}).\leqno(3.5)
$$
Then there exists a $\tau \in (-r^2,-\a r^2)$, such that
$$
\mu(\tau)\geq (\fr 12-\a)(1-\a)^{-1}
mes(S_{r})mes(K_{r}).\leqno(3.6)
$$
From (3.2) and (3.6), we have by noticing $v=0$ when $u\geq 1$
$$
\int_{K_{r}}\int_{S_{\beta r}} v(\tau,x',\overline x)d \overline x
dx'\leq \fr 12(1-\a)^{-1}mes(S_{r})mes(K_{r})\ln(\fr{1}{h^{\fr
98}}).\leqno(3.7)
$$
Now we choose $\a$ (near zero) and $\beta$ (near one), so that
$$
\fr{1}{4\beta^{2Q}}+\fr{1}{2\beta ^{2Q}(1-\a)}\leq \fr 4
5.\leqno(3.8)
$$
By (3.2),( 3.4), (3.7) and (3.8), we deduce
$$\ba{lll}
\ln(\fr 1 {2h})mes(K_{\beta r}\times S_{\beta r}\setminus {\cal
N}_t)\leq &[C(1-\beta)^{-2}\beta^{-Q}+\cr\cr &+\fr 45\ln(\fr
{1}{h^{\fr 98}})]mes(K_{\beta r}\times S_{\beta r}).\ea\leqno(3.9)
$$
Since
$$
\fr{\ln(h^{-\fr 98})}{\ln(h^{-1})}\lra \fr 98,\qquad\hbox{as}
\quad h\ra 0,
$$
then there exists constant $h_1$ such that for $0<h<h_1$ and $t
\in[-\a r^2,0]$
$$
mes(K_{\beta r}\times S_{\beta r}\setminus {\cal N}_t)\leq \fr
{10}{11}mes(K_{\beta r}\times S_{\beta r}).
$$
Then we proved our lemma.
\begin{cor}
Under the assumptions of Lemma 3.1, we can choose $\theta$,
$0<\theta< \a$ and $\theta< \beta$ small enough so that
$$
mes\{{\cal B}^-_{\beta r} \setminus {\cal B}^-_{\theta r} \cap
\{(t,x)|\quad u\geq h\}\}\geq C_0(\a,\beta,\Lambda) mes \{{\cal
B}^-_{\beta r}\},
$$
where $0<C_0(\a,\beta,\Lambda)<1$.
\end{cor}

Let $\chi(s)$ be a smooth function given by
$$\ba{ll}
\chi(s)=1 \qquad if \quad s\leq {\sqrt\theta} r,\\
\chi(s)=0 \qquad if \quad s>\beta r, \ea
$$
where ${\sqrt\theta}<\fr {\beta}{2}$ is a constant. Moreover, we
assume that
$$
0\leq -\chi'(s) \leq \fr{2}{(\beta -{\sqrt\theta})r},
$$
and $\chi'(s)<0$ if ${\sqrt\theta} r<s<\beta r$. We set
$$
\phi_0(t,x)=\chi([|t|^Q\langle
C^{-1}(|t|)e^{tB^T}x,e^{tB^T}x\rangle+\sum_{i=m_0}^{N} \fr
{x_i^2}{r^{2\a_i-2Q}}-c_1 t r^{2Q-2} ]^{\fr {1} {2Q}}),
$$
$$
\phi_1(x)=\chi(\theta |x'|),
$$
$$
\phi(t,x)=\phi_0(t,x)\phi_1(x),\leqno(3.10)
$$
where $c_1>1$ is chosen so that
$$
|2\sum x_ib_{ij} \fr {x_j^2}{r^{2\a_j-2Q}}|+|t|^Q\langle
A_0e^{tB}C^{-1}e^{tB^T}x,e^{tB}C^{-1}e^{tB^T}x\rangle< c_1r^{2Q-2},
$$
for $-r^2 \leq t\leq 0$ and $x \in K_r \times S_r$.

We now have the following Poincare's type inequality.
\begin{lemma}
Let $w$ be a non-negative weak sub-solution of (1.2) in ${\cal
B}_1$. Then there exists a constant $C$, only depends on $\la$ and
$N$, such that for $r<\theta<1$
$$
\int_{{\cal B}^-_{\theta r}}(w(z)-I_0)_+^2\leq
C\fr{r^2}{(1-\theta)^2}\int_{{\cal B}^-_{\fr r
{\theta}}}|D_{m_0}w|^2, \leqno(3.11)
$$
where $I_0$ is given by
$$
I_0=max_{{\cal B}^-_{\theta r}}[I_1(z)+C_2(z)],\leqno(3.12)
$$
and
$$
I_1(z)=\int_{{\cal B}^-_{\fr r {\theta}}} [\langle
{\phi}_1A_0D{\phi}_0,D\Gamma_0(z,\cdot)\rangle
w-\Gamma_0(z,\cdot)wY\phi](\zeta)d\zeta,\leqno(3.13)
$$
$$
C_2(z)=\int_{{\cal B}^-_{\fr r {\theta}}} [\langle
{\phi}_0A_0D{\phi}_1,D\Gamma_0(z,\cdot)\rangle w](\zeta)d\zeta,
$$
where $\Gamma_0$ is the fundamental solution,  and $\phi$ is given
by (3.10).
\end{lemma}
 {\it Proof:} We represent $w$ in terms of the fundamental
solution of $\Gamma_0$. For $z \in {\cal B}^-_{\theta r}$, we have
$$\ba{llll}
w(z)&=\int_{{\cal B}^-_{\fr r {\theta}}}  [\langle
A_0D(w\phi),D\Gamma_0(z,\cdot)\rangle
-\Gamma_0(z,\cdot)Y(w\phi)](\zeta)d\zeta \\ \\&=
I_1(z)+I_2(z)+I_3(z)+C_2(z),\ea\leqno(3.14)
$$
where $I_1(z)$ is given by (3.13) and
$$
I_2(z)=\int_{{\cal B}^-_{\fr r {\theta}}} [\langle
(A_0-A)Dw,D\Gamma_0(z,\cdot)\rangle\phi-\Gamma_0(z,\cdot)\langle
ADw,D\phi\rangle](\zeta)d\zeta,
$$
$$
I_3(z)=\int_{{\cal B}^-_{\fr r {\theta}}} [\langle
ADw,D(\Gamma_0(z,\cdot)\phi)\rangle-\Gamma_0(z,\cdot)\phi
Yw](\zeta)d\zeta,
$$
And $C_2(z)$ denotes the remaining parts of the integral of $I_1$.

From our assumption that $w$ is a weak sub-solution of (1.2), then
$I_3(z)\leq 0$ (see[11]). Then in ${\cal B}^-_{\theta r}$,
$$
0\leq (w(z)-I_0)_+\leq I_2(z).
$$
By Lemma 2.2 we have
$$
||I_2||_{L^2({\cal B}^-_{\theta r})}\leq \theta r||I_2||_{L^{2+\fr 4
Q}({\cal B}^-_{\theta r})}\leq \fr{C
\theta^2r}{1-\theta}||D_{m_0}w||_{L^2({\cal B}^-_{{\fr r
{\theta}}})}.\leqno(3.15)
$$
Then we proved our lemma.

Now we apply Lemma 3.2. to the function
$$
w= ln^+\fr{h}{u+h^{\fr 98}}.
$$
We estimate the value of $I_0$ given by (3.12) and (3.13) in Lemma
3.2.
\begin{lemma}
Under the assumptions of Lemma 3.2, there exist constants
$\lambda_0$, $r_0$ and $r_0<\theta$ only depend on constants $\a$,
 $\beta$, $\Lambda$, $\lambda$, $\phi$ and $\theta$, $0<\lambda_0<1$, such
 that for $r<r_0$
$$
|I_0|\leq \lambda_0 ln(\fr {1}{h^{\fr 18}}).\leqno(3.16)
$$
\end{lemma}
{\it Proof:} We first note that the support of the function $\phi_0
D\phi_1$ is contained in the set ${\cal B}^-_{r}\cap \{(t,x',\ov
x)|\quad{\fr r {\sqrt\theta}}<|x'|<{\fr r {\theta}}\}$. Then it is
easy to check that there exists a constant $C$ which only depends on
$\lambda$ and $N$ such that, for $z \in {\cal B}^-_{\theta r}$
$$
|C_2(z)|\leq C \theta ln(\fr {1}{h^{\fr 18}}).
$$
Therefore $C_2(z)\ra 0$ as $\theta \ra 0$.

Now we let $w\equiv 1$ then (3.14) gives, for $z \in {\cal
B}^-_{\theta r}$,
$$
1=\int_{{\cal B}^-_{\fr r {\theta}}} [\langle
\phi_1A_0D\phi_0,D\Gamma_0(z,\cdot)\rangle
-\Gamma_0(z,\cdot)Y\phi](\zeta)d\zeta+C_2(z)|_{w=1},\leqno(3.17)
$$
where $\phi$ is given by (3.10). Since the matrix $C^{-1}(t)$ is
positive definite for $t>0$ and by the assumption of matrix $B$,
one can check
$$
Y\langle C^{-1}(|t|)x,x\rangle=-\langle
A_0C^{-1}(|t|)x,C^{-1}(|t|)x\rangle,
$$
$$
Y\langle C^{-1}(|t|)e^{tB^T}x,e^{tB^T}x\rangle=-\langle
A_0e^{tB}C^{-1}(|t|)e^{tB^T}x,e^{tB}C^{-1}(|t|)e^{tB^T}x\rangle,\leqno(3.18)
$$
then by the choosing of $c_1$, it is easy to see that
$$
Y\phi \leq 0.\leqno(3.19)
$$

For $z=0$, by our construction of $\phi$, we have
$$
\langle \phi_1 A_0D\phi_0,D\Gamma_0(z,\cdot)\rangle\geq 0,
$$
therefore
$$
\langle \phi_1A_0D\phi_0,D\Gamma_0(z,\cdot)\rangle
-\Gamma_0(z,\cdot)Y\phi \geq 0.\leqno(3.20)
$$
In fact,
$$
\langle A_0D\langle C^{-1}(|t|)e^{tB^T}x,e^{tB^T}x\rangle,D\langle
C^{-1}(|t|)e^{tB^T}x,e^{tB^T}x\rangle\rangle  \geq 0.
$$
We note that the support of $\chi'(s)$ is in the region
${\sqrt\theta}r<s<\beta r$. Thus for some ${\beta}'<\beta$, the set
${\cal B}^-_{{\beta}' r} \setminus {\cal B}^-_{{\sqrt\theta} r}$
with $|t|>\theta^2 r^2$ is contained in the support of $\phi'$ and
then the inequality holds in (3.20). By the choosing of $c_1$, we
know that (3.20) is positive in ${\cal B}^-_{{\beta}' r} \setminus
{\cal B}^-_{{\sqrt\theta} r}$ with $|t|>\theta^2 r^2$. Then the
integral of (3.20) on the domain ${\cal B}^-_{{\beta}' r} \setminus
{\cal B}^-_{{\sqrt\theta} r}$ with $|t|>\theta^2 r^2$ is lower
bounded by a positive constant which independent of small $r$ and
$\theta$.

Then we may choose $\theta$ small so that (3.20) still holds for $z
\in {\cal B}^-_{\theta r}$ and the inequality holds in ${\cal
B}^-_{{\beta}' r} \setminus {\cal B}^-_{{\sqrt\theta} r}$ with
$|t|>\theta^2 r^2$. Therefore the integral function in (3.17) is
nonnegative for $z \in {\cal B}^-_{\theta r}$ and positive in ${\cal
B}^-_{{\beta}' r} \setminus {\cal B}^-_{{\sqrt\theta} r}$ with
$|t|>\theta^2 r^2$. Since $w=0$ when $u\geq h$ and $w\leq ln(\fr
{1}{h^{\fr 18}})$, then by choosing $\theta$ small enough, our lemma
follows from Corollary 3.1 and (3.17).

\begin{lemma}
Suppose that $u(t,x)\geq 0$ be a solution of equation (1.2) in
${\cal B}^-_r$ centered at $(0,0)$ and
$$
mes\{(t,x)\in {\cal B}^-_r, \quad u \geq 1\} \geq \fr 1 2 mes
({\cal B}^-_r).
$$
Then there exist constant $\theta$ and $h_0$, $0<\theta, h_0<1$
which only depend on $\la$, $\la_0$ and $N$ such that
$$
u(t,x) \geq h_0\quad \hbox{in}\quad {\cal B}^-_{\theta r}.
$$
\end{lemma}
{\it Proof:} We consider $$w=\ln^+(\fr{h}{u+h^{\fr98}}),$$ for
$t\in[-\a r^2,0]$. By applying Lemma 3.2 to $w$ and as in the
proof of Lemma 3.1, we have
$$
-\!\!\!\!\!\!\int_{{\cal B}^-_{\theta r}}( w-I_0)_+^2 \leq C
(ln(h^{-\fr 18})) .\leqno(3.21)
$$
By Lemma 2.1, there exists a constant, still denoted by $\theta$,
such that for $z \in {\cal B}^-_{\theta r}$,
$$
w-I_0\leq C (ln(h^{-\fr 18}))^{\fr 12} .\leqno(3.22)
$$
Therefore we may choose $h_0$ small enough, so that
$$
C (\ln (\fr {1}{h_0^{\fr 18}}))^{\fr 12}\leq  \ln (\fr
{1}{2h_0^{\fr 18}})-\lambda_0\ln (\fr {1}{h_0^{\fr 18}}).
$$
Then (3.16) and (3.22) implies
$$
\max_{{\cal B}^-_{\theta r}}\fr{h_0}{u+h_0^{\fr 98}}\leq \fr
{1}{2h_0^{\fr 18}},
$$
which implies $\min_{{\cal B}^-_{\theta r}}u\geq h_0^{\fr 98}$.
Then we finished the proof of Lemma.

{\bf Proof of Theorem 1.1.} We may assume that $M=\max_{{\cal
B}^-_{r}}\pm u$, otherwise we replace $u$ by $u-c$. Then either
$1+\fr u M$ or $1-\fr u M$ satisfies the assumption of Lemma 3.4,
thus Lemma 3.4 implies
$$
Osc_{{\cal B}^-_{\theta r}}u\leq (1-\fr{h_0}{2})Osc_{{\cal
B}^-_{r}}u,
$$
which implies the $C^{\a}$ regularity of $u$ near point $(0,0)$ by
the standard iteration arguments. By the left invariant
translation group action, we know that $u$ is $C^{\a}$ in the
interior.

{\bf Acknowledgements}: The research is partially supported by the
Chinese NSF under grant 10325104 and the innovation program at
CAS. The author thanks Prof. Xin Zhouping for many valuable
discussions on this subject.

\end{document}